\documentclass[11pt]{amsart}

\makeatletter
\@namedef{subjclassname@2020}{%
  \textup{2020} Mathematics Subject Classification}
\makeatother

\usepackage{amsmath,amsfonts,amssymb,amsthm,bm}

\title[Galois lines for the maximal curve with the second largest genus]{Galois lines for the maximal curve with the second largest genus in characteristic two}
\author{Kazuki Higashine}

\address[Kazuki Higashine]{Miyako Junior College, Iwate Prefectural University, Kanan 1-5-1, Miyako 027-0039, Japan}
\email{k.higashine.1029@gmail.com}

\newtheorem{theorem}{Theorem}

\newtheorem{lemma}[theorem]{Lemma}
\newtheorem{fact}[theorem]{Fact}
\newtheorem*{mtheorem}{Theorem}
\newtheorem*{coro}{Corollary}
\newtheorem*{mondai}{Problem}

\theoremstyle{definition}

\begin{document}
\begin{abstract}
We determine all Galois lines for the embedding of the maximal curve with the second largest genus in characteristic two into $\mathbb{P}^{3}$ induced by its Frobenius linear series. We also determine the associated Galois groups.
\end{abstract}
\maketitle

\section{Introduction}
Let $\mathcal{X}$ be a non-singular, projective, irreducible algebraic curve of genus $g=g(\mathcal{X})$ defined over an algebraically closed field $k$ of characteristic $p\geq0$. Let $\varphi :\mathcal{X}\to \mathbb{P}^{2}$ be a birational embedding, that is, a morphism which is birational onto its image $\varphi (\mathcal{X})$. A point $P\in \mathbb{P}^{2}$ is called a Galois point (for $\varphi (\mathcal{X})$) if the field extension $k(\varphi (\mathcal{X}))/\pi_{P}^{\ast}k(\mathbb{P}^{1})$ of function fields induced by the projection $\pi_{P}:\varphi (\mathcal{X})\dashrightarrow \mathbb{P}^{1}$ from $P$ is a Galois extension. The notion of a Galois point was introduced by Yoshihara in 1996 \cite{FukasawaDedicata2009, MiuraYoshiharaAlgebra2000, YoshiharaAlgebra2001}. Analogously, the notion of a Galois line was introduced by Yoshihara. Let $\varphi :\mathcal{X}\to \mathbb{P}^{3}$ be a birational embedding. A line $\ell \subset \mathbb{P}^{3}$ is called a Galois line if the field extension $k(\varphi (\mathcal{X}))/\pi_{\ell}^{\ast}k(\mathbb{P}^{1})$ of function fields induced by the projection $\pi_{\ell}:\varphi (\mathcal{X})\dashrightarrow \mathbb{P}^{1}$ from $\ell$ is a Galois extension \cite{DuyaguitYoshiharaAlgebraColloq.2005, YoshiharaAlgebraColloq.2006}. A Galois line $\ell$ is called a skew Galois line if $\ell \cap \varphi (\mathcal{X})=\emptyset$. For a Galois line $\ell$, the associated Galois group is denoted by $G_{\ell}$. Let $q$ be a power of $p>0$. A curve $\mathcal{X}$ defined over $\mathbb{F}_{q^{2}}$ is called $\mathbb{F}_{q^{2}}$-maximal if $|\mathcal{X}(\mathbb{F}_{q^{2}})|=q^{2}+1+2gq$. The Hermitian curve $\mathcal{H}_{q}$ defined over the finite field $\mathbb{F}_{q^{2}}$ is a classical example of a maximal curve. For $q>2$, it is, up to projective equivalence, the unique non-singular plane curve of degree $q+1$ attaining the maximum possible number of Galois points \cite{FukasawaPadova2013}. In this case, the set of all $\mathbb{F}_{q^{2}}$-rational points of $\mathbb{P}^{2}$ coincides with the set of all Galois points for the Hermitian curve $\mathcal{H}_{q}$ \cite{HommaComm.2006}. It is also well known that the Hermitian curve $\mathcal{H}_{q}$ attains the maximum possible genus among all $\mathbb{F}_{q^{2}}$-maximal curves \cite{Ihara1981}. By the natural embedding theorem for a maximal curve $\mathcal{X}$ defined over $\mathbb{F}_{q^2}$ (see \cite{FuhrmannTorresRend.Circ.Mat.PalermoSuppl.1998}, \cite[Section 10.3]{hkt}, \cite{KorchmarosTorresCompositioMath.2001}), the Hermitian curve is the only maximal curve whose Frobenius linear series $|(q+1)P|$, where $P$ is an $\mathbb{F}_{q^{2}}$-rational point of the curve, defines an embedding into $\mathbb{P}^{2}$ (see \cite[Theorem 2.4]{FuhrmannTorresRend.Circ.Mat.PalermoSuppl.1998}, \cite{RuckStichtenothJ.ReineAngew.Math.1994}). If $\mathcal{X}$ is an $\mathbb{F}_{q^{2}}$-maximal curve with $g(\mathcal{X})<g(\mathcal{H}_{q})$, then $g(\mathcal{X})\leq \lfloor (q-1)^{2}/4\rfloor$ and $\mathcal{X}$ can be embedded into $\mathbb{P}^{r}$ for some $r\geq 3$ by the Frobenius linear series (see \cite{FuhrmannTorresManuscriptaMath.1996}, \cite[Theorem 2.4]{FuhrmannTorresRend.Circ.Mat.PalermoSuppl.1998}, \cite[Corollary 10.25]{hkt}). In particular, any $\mathbb{F}_{q^{2}}$-maximal curve $\mathcal{X}$ of genus $g(\mathcal{X})=\lfloor (q-1)^{2}/4\rfloor$ is birationally equivalent to one of the following plane curves:
\begin{enumerate}
\item[(1)] $y^{\frac{q+1}{2}}=x^{q}+x$, for odd $q$,
\item[(2)] $y^{q+1}=x+x^{2}+\cdots +x^{\frac{q}{2}}$, for even $q$
\end{enumerate}
(see \cite{AbdonTorresManuscriptaMath.1999, FuhrmannGarciaTorresJ.NumberTheory1997, hkt, KorchmarosTorresMath.Ann.2002}, \cite[p.244]{MMR2008}). The non-singular projective models of curves (1) and (2) can be embedded into $\mathbb{P}^{3}$ via their Frobenius linear series. These observations lead us to the following general problem.
\begin{mondai}
Let $\mathcal{X}$ be a non-singular projective $\mathbb{F}_{q^{2}}$-maximal curve. Suppose that its Frobenius linear series induces an embedding $\varphi:\mathcal{X}\to\mathbb{P}^{3}$. Describe the arrangement of all Galois lines for $\varphi(\mathcal{X})$.
\end{mondai}
\noindent We focus on $\mathbb{F}_{q^{2}}$-maximal curves with the second largest possible genus $\lfloor (q-1)^{2}/4\rfloor$. For odd $q\geq 5$, Fukasawa determined the arrangement of all Galois lines for the curve (1) \cite{FukasawaFFTA2026}. In this article, we determine all Galois lines for the curve (2), where $q=2^n\geq 4$.

From now on, we assume that $p=2$ and $q=2^{n}\geq 4$. We denote by $\mathcal{X}$ the non-singular projective model of the affine curve (2). We naturally identify the function field $k(\mathcal{X})$ with the function field $k(x, y)$ of the affine curve (2). Let $P_{\infty}\in \mathcal{X}(\mathbb{F}_{q^{2}})$ be the unique pole of both $x$ and $y$. We consider the morphism $\varphi : \mathcal{X}\to \mathbb{P}^{3}$ defined by the complete linear series $|(q+1)P_{\infty}|$, given by $\varphi=(x:y:1:y^2)$. Our main result is the following theorem.

\begin{mtheorem}\label{maintheorem}
Let $(X:Y:Z:W)$ be a system of homogeneous coordinates on $\mathbb{P}^{3}$, and let $\ell \subset \mathbb{P}^{3}$ be a line. Then $\ell$ is a Galois line for $\varphi (\mathcal{X})$ if and only if $\ell$ is one of the following:
\begin{enumerate}
\item[(1)] the line $\ell_{b}^{(1)}$ defined by $X+b^{q}Y+b^{2q}W=Z=0$ for some $b\in \mathbb{F}_{q^{2}}$,
\item[(2)] the line $\ell_{c}^{(2)}$ defined by $Y+cZ=W+c^{2}Z=0$ for some $c\in k$,
\item[(3)] the line $\ell_{d}^{(3)}$ defined by $Y+dW=Z=0$ for some $d\in \mathbb{F}_{q^{2}}$.
\end{enumerate}
In case $({\rm 1})$, $G_{\ell_{b}^{(1)}}\simeq C_{q+1}$, where $C_{q+1}$ is a cyclic group of order $q+1$. In case $({\rm 2})$, $G_{\ell_{c}^{(2)}}\simeq (\mathbb{Z}/2\mathbb{Z})^{\oplus n-1}$. In case $({\rm 3})$, the Galois group is determined as follows:
\begin{enumerate}
\item[$\bullet$] $G_{\ell_{0}^{(3)}}\simeq (\mathbb{Z}/2\mathbb{Z})^{\oplus n-1}$,
\item[$\bullet$] $G_{\ell_{d}^{(3)}}\simeq (\mathbb{Z}/2\mathbb{Z})^{\oplus n}$ if $d\not= 0$ and $d^{q+1}=1$,
\item[$\bullet$] $G_{\ell_{d}^{(3)}}\simeq (\mathbb{Z}/2\mathbb{Z})^{\oplus n-2}\oplus \mathbb{Z}/4\mathbb{Z}$ if $d\not= 0$ and $d^{q+1}\not= 1$.
\end{enumerate}
\end{mtheorem}

\section{Preliminaries}
We recall some properties of $\varphi (\mathcal{X})$. It is easily verified that
\[ (x)_{\infty}=(q+1)P_{\infty},\quad (y)_{\infty}=\frac{q}{2}P_{\infty}. \]
The genus $g(\mathcal{X})$ is equal to $\left\lfloor (q-1)^{2}/4\right\rfloor=(q(q-2))/4$ (see \cite{AbdonTorresManuscriptaMath.1999}, \cite[Lemma 12.1]{hkt}, \cite{StichtenothArch.1973}). Let $P_{0}\in \mathcal{X}(\mathbb{F}_{q^{2}})$ be the common zero of $x$ and $y$. Then $\varphi (P_{0})=(0:0:1:0)$ and $\varphi (P_{\infty})=(1:0:0:0)$.

For a hyperplane $H\subset \mathbb{P}^{3}$, the intersection divisor of $\varphi (\mathcal{X})$ and $H$ on $\mathcal{X}$ is denoted by $\varphi^{\ast}H$. For a point $P\in \mathcal{X}$, the order of $\varphi^{\ast}H$ at $P$ is denoted by ${\rm ord}_{P}(\varphi^{\ast}H)$. For each $P\in \mathcal{X}$, there exists a unique hyperplane $H_{P}\subset \mathbb{P}^{3}$ with the maximum order at $P$. The hyperplane $H_{P}$ is called the osculating hyperplane at $P$. The osculating hyperplane at $P\in \mathcal{X}\setminus \{ P_{\infty}\}$ is defined by $X+b^{q}Y+a^{q}Z+b^{2q}W=0$ if $\varphi (P)=(a:b:1:b^2)$. In fact, the equality $x+x^{q}=y^{q+1}+y^{2q+2}$ holds as an element of $k(\mathcal{X})=k(x, y)$. Since
\[ \begin{aligned}
x+b^{q}y+a^{q}+b^{2q}y^{2}&=x^{q}+y^{q+1}+y^{2q+2}+b^{q}y+a^{q}+b^{2q}y^{2}\\
&=x^{q}+a^{q}+y^{q+1}+b^{q}y+y^{2q+2}+b^{2q}y^{2}\\
&=(x+a)^{q}+y(y+b)^{q}+y^{2}(y+b)^{2q}
\end{aligned}\]
holds as an element of $k(x, y)$, it follows ${\rm ord}_{P}(\varphi^{\ast}H_{P})\geq q$. If $P\in \mathcal{X}\setminus \mathcal{X}(\mathbb{F}_{q^{2}})$, ${\rm ord}_{P}(\varphi^{\ast}H_{P})=q$ holds since the point $\varphi (P^{q^{2}})=(a^{q^{2}}:b^{q^{2}}:1:b^{2q^{2}})$ is contained in the hyperplane $H_{P}$. For each $P\in \mathcal{X}(\mathbb{F}_{q^{2}})\setminus \{ P_{\infty}\}$, it follows ${\rm ord}_{P}(\varphi^{\ast}H_{P})=q+1$ (see \cite[Section 10.3]{hkt}). The osculating hyperplane at $P_{\infty}$ is defined by $Z=0$. In particular, $H_{P_{0}}$ is defined by $X=0$. Since the order of the hyperplane defined by $W=0$ at $P_{0}$ is equal to $2$, it follows that the (projective) tangent line $T_{P_{0}}$ at $P_{0}$ is defined by $X=W=0$. Similarly, since the order of the hyperplane defined by $Y=0$ at $P_{\infty}$ is equal to $q/2+1$, it follows that the tangent line $T_{P_{\infty}}$ at $P_{\infty}$ is defined by $Y=Z=0$. We summarize the hyperplane orders of $\varphi (\mathcal{X})$.

\begin{fact}[\cite{AbdonTorresManuscriptaMath.1999}, Proposition 2.4, \cite{hkt}, Section 10.4]\label{orderseq}
Let $P\in \mathcal{X}$ and let $H\subset \mathbb{P}^{3}$ be a hyperplane.
\begin{enumerate}
\item[(1)] If $P=P_{\infty}$, then ${\rm ord}_{P}(\varphi^{\ast}H)=0, 1, q/2+1$ or $q+1$.
\item[(2)] If $P\in \mathcal{X}(\mathbb{F}_{q^{2}})\setminus \{ P_{\infty}\}$, then ${\rm ord}_{P}(\varphi^{\ast}H)=0, 1, 2$ or $q+1$.
\item[(3)] If $P\in \mathcal{X}\setminus \mathcal{X}(\mathbb{F}_{q^{2}})$, then ${\rm ord}_{P}(\varphi^{\ast}H)=0, 1, 2$ or $q$.
\end{enumerate}
\end{fact}

Let $f\in \mathcal{L}((q+1)P_{\infty})$, where $\mathcal{L}((q+1)P_{\infty})$ is the $k$-vector space associated to the divisor $(q+1)P_{\infty}$. Since $f$ is regular at $\mathcal{X}\setminus \{ P_{\infty}\}$, $f$ is a polynomial in $x$ and $y$. By comparing the orders at $P_{\infty}$ of the monomials appearing in $f$, we have the following fact.

\begin{fact}[\cite{hkt}, Lemma 12.2, \cite{StichtenothArch.1973}]\label{complete}
We have $\mathcal{L}((q+1)P_{\infty})=\left\langle x, y, 1, y^{2}\right\rangle$.
\end{fact}

The canonical divisor $K$ of $\mathcal{X}$ is given by $((q^{2}-2q-4)/2)P_{\infty}$ (see \cite[Lemma 12.1]{hkt}, \cite{StichtenothArch.1973}). For $P\in \mathcal{X}\setminus \{ P_{\infty}\}$ which satisfies $\varphi (P)=(a:b:1:b^{2})$, the function $f=(y+b)^{\frac{q}{2}-1}$ is contained in $\mathcal{L}(K)$. The order of $f$ at $P$ is equal to $q/2-1$. By the relation between the gap sequence and the vanishing orders of holomorphic differentials, $q/2$ is a gap at every point $P\not= P_{\infty}$. Therefore, we have the following.

\begin{fact}[\cite{hkt}, Theorem 12.11, \cite{StichtenothArch.1973}]\label{fixinfinite}
For any $\sigma \in {\rm Aut}(\mathcal{X})$, $\sigma (P_{\infty})=P_{\infty}$.
\end{fact}

We identify ${\rm Aut}(\mathbb{P}^{3})$ with the projective linear group ${\rm PGL}(4, k)$. Let $\sigma \in {\rm Aut}(\mathcal{X})$. By Facts \ref{complete} and \ref{fixinfinite}, the linear system corresponding to the morphism $\varphi \circ \sigma$ is $\mathcal{D}_{\mathcal{X}}=|(q+1)P_{\infty}|$. Therefore, there exists an injective homomorphism of groups
\[ {\rm Aut}(\mathcal{X}) \to {\rm PGL}(4, k), \]
of which the image preserves the curve $\varphi (\mathcal{X})\subset \mathbb{P}^{3}$. Hereafter, we identify ${\rm Aut}(\mathcal{X})$ with the subgroup of ${\rm PGL}(4, k)$ by the above homomorphism. 

Put $L_{q}(T)=T+T^{2}+\cdots +T^{\frac{q}{4}}+T^{\frac{q}{2}}$. We note that, for any $b\in \mathbb{F}_{q^{2}}$, the equation $b^{q+1}=L_{q}(T)$ has exactly $q/2$ distinct solutions in $\mathbb{F}_{q^{2}}$. Indeed, since $L_{q}(T)^{2}+L_{q}(T)=T^{q}+T$, any solution $a\in k$ of $b^{q+1}=L_{q}(a)$ satisfies $b^{2q+2}+b^{q+1}=a^{q}+a$. Since $(b^{2q+2}+b^{q+1})^{q}=b^{2q+2}+b^{q+1}$, it follows $a^{q^{2}}=a$. Since $L_{q}(T)-b^{q+1}$ is a separable polynomial of degree $q/2$, it has exactly $q/2$ distinct roots in $k$, all of which belong to $\mathbb{F}_{q^{2}}$. These computations show that the set $\varphi (\mathcal{X}(\mathbb{F}_{q^{2}}))\setminus \{ (1:0:0:0)\}$ coincides with the set $\{ (a:b:1:b^{2})\in \mathbb{P}^{3}\mid a, b\in \mathbb{F}_{q^{2}}, b^{q+1}=L_{q}(a)\}$. We consider the following subgroup $N$ of ${\rm PGL}(4, k)$:
\[ N=\left\{ \left(
\begin{array}{cccc}
1 & b^{q} & a & b^{2q} \\
0 & 1 & b & 0 \\
0 & 0 & 1 & 0 \\
0 & 0 & b^{2} & 1
\end{array}
\right) \ \middle| \ a, b\in \mathbb{F}_{q^{2}}, b^{q+1}=L_{q}(a)\right\}. \]
A direct computation shows that $N$ preserves $\varphi (\mathcal{X})$ and has order $q^{3}/2$. For any point $(a:b:1:b^{2})\in \varphi (\mathcal{X}(\mathbb{F}_{q^{2}}))\setminus \{ (1:0:0:0)\}$,
\[ \left(
\begin{array}{cccc}
1 & b^{q} & a & b^{2q} \\
0 & 1 & b & 0 \\
0 & 0 & 1 & 0 \\
0 & 0 & b^{2} & 1
\end{array}
\right) \left(
\begin{array}{cccc}
0 \\
0 \\
1 \\
0
\end{array}
\right) =\left(
\begin{array}{cccc}
a \\
b \\
1 \\
b^{2}
\end{array}
\right) \]
holds. Since linear automorphisms preserve hyperplane orders, Fact \ref{orderseq} implies that ${\rm Aut}(\mathcal{X})$ preserves the set $\mathcal{X}(\mathbb{F}_{q^{2}})\setminus \{ P_{\infty}\}$. Therefore, the following holds.

\begin{fact}\label{transitive}
${\rm Aut}(\mathcal{X})$ acts on the set $\mathcal{X}(\mathbb{F}_{q^{2}})\setminus \{ P_{\infty}\}$ transitively.
\end{fact}

For a finite subgroup $G\subset {\rm Aut}(\mathcal{X})$ and a point $P\in \mathcal{X}$, the stabilizer of $P$ in $G$ (resp. the orbit of $P$ under $G$) is denoted by $G(P)$ (resp. $G\cdot P$). By Fact \ref{transitive} and the Orbit-Stabilizer theorem,
\[ [{\rm Aut}(\mathcal{X}) : {\rm Aut}(\mathcal{X})(P_{0})]=|{\rm Aut}(\mathcal{X})\cdot P_{0}|=|\mathcal{X}(\mathbb{F}_{q^{2}})\setminus \{ P_{\infty}\}|=\frac{q^{3}}{2}. \]

\begin{fact}\label{kernel}
The stabilizer ${\rm Aut}(\mathcal{X})(P_{0})$ is given by
\[ H=\left\{ \left(
\begin{array}{cccc}
1 & 0 & 0 & 0 \\
0 & \zeta & 0 & 0 \\
0 & 0 & 1 & 0 \\
0 & 0 & 0 & \zeta^{2}
\end{array}
\right) \ \middle| \ \zeta \in k, \zeta^{q+1}=1 \right\} . \]
In particular, any element of $H$ is defined over $\mathbb{F}_{q^{2}}$.
\end{fact}

\begin{proof}
Let $\tau \in {\rm Aut}(\mathcal{X})(P_{0})$. Then $\tau (P_{0})=P_{0}$ holds. By Fact \ref{fixinfinite}, $\tau (P_{\infty})=P_{\infty}$ holds. As noted above, the osculating hyperplane $H_{P_{0}}$ (resp. $H_{P_{\infty}}$) at $P_{0}$ (resp. $P_{\infty}$) is defined by $X=0$ (resp. $Z=0$). The (projective) tangent line $T_{P_{0}}$ (resp. $T_{P_{\infty}}$) at $P_{0}$ (resp. $P_{\infty}$) is defined by $X=W=0$ (resp. $Y=Z=0$). Since ${\rm Aut}(\mathcal{X})\subset {\rm PGL}(4, k)$, $\tau (T_{P_{0}})=T_{P_{0}}$, $\tau (H_{P_{0}})=H_{P_{0}}$, $\tau (T_{P_{\infty}})=T_{P_{\infty}}$ and $\tau (H_{P_{\infty}})=H_{P_{\infty}}$ hold. Therefore, $\tau$ is represented by the matrix
\[ \tau=\left(
\begin{array}{cccc}
a_{11} & 0 & 0 & 0 \\
0 & a_{22} & 0 & 0 \\
0 & 0 & 1 & 0 \\
0 & 0 & 0 & a_{44}
\end{array}
\right) , \]
where $a_{11}, a_{22}, a_{44}\in k\setminus \{ 0\}$. Since $\tau$ acts on $\varphi (\mathcal{X})$,
\[ a_{22}^{2}=a_{44},\quad a_{11}x+a_{11}^{2}x^{2}+\cdots +a_{11}^{\frac{q}{2}}x^{\frac{q}{2}}=a_{22}^{q+1}y^{q+1}. \]
The relation $y^{q+1}=x+x^{2}+\cdots +x^{\frac{q}{2}}$ implies that
\[ a_{11}x+a_{11}^{2}x^{2}+\cdots +a_{11}^{\frac{q}{2}}x^{\frac{q}{2}}=a_{22}^{q+1}x+a_{22}^{q+1}x^{2}+\cdots +a_{22}^{q+1}x^{\frac{q}{2}}. \]
Since $x$ is algebraically independent over $k$, $a_{11}=a_{11}^{2}=a_{22}^{q+1}$ holds. Therefore, $a_{11}=1$ and $a_{22}^{q+1}=1$ hold. Since $a_{22}^{2}=a_{44}$, $\tau \in H$. $H\subset {\rm Aut}(\mathcal{X})(P_{0})$ is obvious. If $\zeta \in k$ satisfies $\zeta^{q+1}=1$, then $\zeta^{q^{2}}=\zeta^{q^{2}-1}\zeta=(\zeta^{q+1})^{q-1}\zeta=\zeta$. Therefore, every element of $H$ is defined over $\mathbb{F}_{q^{2}}$.
\end{proof}

By Fact \ref{kernel}, $|{\rm Aut}(\mathcal{X})|=(q^{3}(q+1))/2$. Since $H$ normalizes $N$ and $N\cap H=\{ 1\}$, the subgroup generated by $N$ and $H$ is their semidirect product $N\rtimes H$. Since $|N\rtimes H|=(q^{3}(q+1))/2$, ${\rm Aut}(\mathcal{X})=N\rtimes H$. Therefore, we have the following.

\begin{fact}\label{auto}
The full automorphism group ${\rm Aut}(\mathcal{X})$ is given by
\[ \left\{ \left(
\begin{array}{cccc}
1 & \zeta b^{q} & a & \zeta^{2}b^{2q} \\
0 & \zeta & b & 0 \\
0 & 0 & 1 & 0 \\
0 & 0 & b^{2} & \zeta^{2}
\end{array}
\right) \ \middle| \ a, b, \zeta \in \mathbb{F}_{q^{2}}, b^{q+1}=L_{q}(a), \zeta^{q+1}=1\right\} . \]
\end{fact}
\noindent The automorphism group of $\mathcal{X}$ is also described in \cite{HuMaFFTA2026}. 

Finally, we recall some properties of a Galois covering for the proof of our main theorem (see \cite[{\rm I\hspace{-.01em}I\hspace{-.01em}I}. 7.1, 7.2 and 8.2]{Stichtenoth}).
\begin{fact}\label{Galoiscov}
Let $\theta :\mathcal{Y}_{1}\to \mathcal{Y}_{2}$ be a surjective morphism of non-singular projective curves, and let the field extension $k(\mathcal{Y}_{1})/\theta^{\ast}k(\mathcal{Y}_{2})$ be a Galois extension with the Galois group $G$. The ramification index of $\theta$ at a point $P\in \mathcal{Y}_{1}$ is denoted by $e_{P}(\theta)$. Then the following hold.
\begin{itemize}
\item[(1)] If $P, Q\in \mathcal{Y}_{1}, \theta (P)=\theta (Q)$, then there exists an element $\sigma \in G$ such that $\sigma (P)=Q$.
\item[(2)] If $P, Q\in \mathcal{Y}_{1}, \theta (P)=\theta (Q)$, then $e_{P}(\theta)=e_{Q}(\theta)$.
\item[(3)] For each point $P\in \mathcal{Y}_{1}$, the order $|G(P)|$ is equal to $e_{P}(\theta)$.
\end{itemize}
\end{fact}

\section{Proof of the main theorem}
First, we prove the following lemma needed later.

\begin{lemma}\label{containinfinite}
If $\ell$ is a non-skew Galois line, then $\ell$ contains the point $\varphi (P_{\infty})=(1:0:0:0)$.
\end{lemma}

\begin{proof}
Let $\ell$ be a non-skew Galois line. Assume by contradiction that $\ell$ does not contain the point $\varphi(P_{\infty})$. By Facts \ref{fixinfinite} and \ref{Galoiscov}, $e_{P_{\infty}}(\pi_{\ell})=|G_{\ell}(P_{\infty})|=|G_{\ell}|={\rm deg}(\pi_{\ell})$. Since $q\geq 4$, we have $g(\mathcal{X})\geq 2$, and hence $1<{\rm deg}(\pi_{\ell})$. Since $\ell$ is non-skew, ${\rm deg}(\pi_{\ell})<q+1$. Let $H$ be the hyperplane spanned by $\ell$ and $\varphi(P_{\infty})$. Since $\varphi (P_{\infty})\not\in \ell$, $e_{P_{\infty}}(\pi_{\ell})={\rm ord}_{P_{\infty}}(\varphi^{\ast}H)=q/2+1$ holds by Fact \ref{orderseq}. Therefore, $q/2+1$ divides $|{\rm Aut}(\mathcal{X})|$. However, since $q/2+1$ is odd and $\gcd(q/2+1, q+1)=1$, $q/2+1$ does not divide $|{\rm Aut}(\mathcal{X})|=(q^{3}(q+1))/2$. This is a contradiction. Therefore, $\ell$ contains the point $\varphi (P_{\infty})$.
\end{proof}

\begin{proof}[{\bf Proof of Theorem}]

First, we consider skew Galois lines for $\varphi (\mathcal{X})$. We show that any line $\ell_{b}^{(1)}$ defined by $X+b^{q}Y+b^{2q}W=Z=0$ for $b\in \mathbb{F}_{q^{2}}$ is a skew Galois line. Let $\ell_{0}^{(1)}$ be the line defined by $X=Z=0$. It is easily verified that the field extension $k(\varphi (\mathcal{X}))/\pi_{\ell_{0}^{(1)}}^{\ast}k(\mathbb{P}^{1})$ coincides with the field extension $k(x, y)/k(x)$. Therefore, $\ell_{0}^{(1)}$ is a skew Galois line with $G_{\ell_{0}^{(1)}}=H\simeq C_{q+1}$. We consider a line $\ell_{b}^{(1)}$ defined by $X+b^{q}Y+b^{2q}W=Z=0$ for $b\in \mathbb{F}_{q^{2}}$, and a point $P\in \mathcal{X}(\mathbb{F}_{q^{2}})\setminus \{ P_{\infty} \}$ which satisfies $\varphi (P)=(a:b:1:b^{2})$ for some $a\in \mathbb{F}_{q^{2}}$. By Fact \ref{transitive}, there exists $\tau \in {\rm Aut}(\mathcal{X})$ such that $\tau (P_{0})=P$, where $\varphi(P_{0})=(0:0:1:0)$. Since ${\rm Aut}(\mathcal{X})\subset {\rm PGL}(4, k)$ and $\tau (P_{\infty})=P_{\infty}$ by Fact \ref{fixinfinite}, it follows that $\tau (H_{P_{0}})=H_{P}$ and $\tau (H_{P_{\infty}})=H_{P_{\infty}}$. Therefore, $\tau (\ell_{0}^{(1)})=\ell_{b}^{(1)}$, and $\ell_{b}^{(1)}$ is a skew Galois line. Conversely, let $\ell$ be a skew Galois line. Since $\ell \cap \varphi (\mathcal{X})=\emptyset$, ${\rm deg}(\pi_{\ell})=q+1$. Since $|G_{\ell}|=q+1$ is odd, $\pi_{\ell}$ is tamely ramified. For any ramification point $P$, the ramification index divides $q+1$ by Fact \ref{Galoiscov}. Since $\gcd(q/2+1, q+1)=1$, it follows from Fact \ref{orderseq} that every ramification point of $\pi_{\ell}$ is $\mathbb{F}_{q^{2}}$-rational and has ramification index $q+1$. Let $r$ be the number of ramification points. By the Riemann-Hurwitz formula, 
\[ \dfrac{q^{2}-2q-4}{2}=2g(\mathcal{X})-2=-2(q+1)+rq.\]
Therefore, $r=q/2+1>1$. Hence, there exists a point $P\in\mathcal{X}(\mathbb{F}_{q^{2}})\setminus\{P_{\infty}\}$ such that $e_P(\pi_{\ell})=q+1$. Therefore, $\ell\subset H_{P}$. By Facts \ref{fixinfinite} and \ref{Galoiscov}, $e_{P_{\infty}}(\pi_{\ell})=q+1$. Therefore, $\ell \subset H_{P_{\infty}}$. Since $H_{P}$ is defined by $X+b^{q}Y+a^{q}Z+b^{2q}W=0$ for some $a, b\in \mathbb{F}_{q^{2}}$, and $H_{P_{\infty}}$ is defined by $Z=0$, $\ell$ is defined by $X+b^{q}Y+b^{2q}W=Z=0$.

Next, we consider non-skew Galois lines for $\varphi (\mathcal{X})$ by dividing into the following two cases $\rm(\,I\,)$ and $\rm(I\hspace{-.01em}I)$. 

$\rm(\,I\,)$ We consider non-skew Galois lines which are not contained in the hyperplane defined by $Z=0$. 
We show that the line $\ell_{c}^{(2)}$ defined by $Y+cZ=W+c^{2}Z=0$ is such a non-skew Galois line for any $c\in k$. It is easily verified that the field extension $k(\varphi (\mathcal{X}))/\pi_{\ell_{c}^{(2)}}^{\ast}k(\mathbb{P}^{1})$ coincides with the field extension $k(x, y)/k((y+c)/(y+c)^{2})=k(x, y)/k(y)$. Therefore, $\ell_{c}^{(2)}$ is a non-skew Galois line with
\[ G_{\ell_{c}^{(2)}}=\left\{ \left(
\begin{array}{cccc}
1 & 0 & \lambda & 0 \\
0 & 1 & 0 & 0 \\
0 & 0 & 1 & 0 \\
0 & 0 & 0 & 1
\end{array}
\right) \ \middle| \ \lambda \in \mathbb{F}_{q^{2}}, L_{q}(\lambda)=0 \right\}\simeq (\mathbb{Z}/2\mathbb{Z})^{\oplus n-1} . \]
Since $(0:c:1:c^{2})\in \ell_{c}^{(2)}$, $\ell_{c}^{(2)}$ is not contained in the hyperplane defined by $Z=0$. Conversely, let $\ell$ be a non-skew Galois line which is not contained in the hyperplane defined by $Z=0$. By Lemma \ref{containinfinite}, $\varphi (P_{\infty})\in \ell$. Since $\ell$ is not contained in the hyperplane defined by $Z=0$, there exist $c,d\in k$ such that $\ell$ is defined by $Y+cZ=W+dZ=0$. The orders at $P_{\infty}$ of these two hyperplane sections are $q/2+1$ and $1$, respectively. It follows that $q/2=(q/2+1)-1=e_{P_{\infty}}(\pi_{\ell})=|G_{\ell}(P_{\infty})|=|G_{\ell}|={\rm deg}(\pi_{\ell})$ by Facts \ref{orderseq}, \ref{fixinfinite} and \ref{Galoiscov}. Since the degree of the fixed divisor of the pencil is $q/2+1$ and its multiplicity at $P_{\infty}$ is $1$, its support contains a point in $\mathcal{X}\setminus \{ P_{\infty} \}$. Therefore, there exists a point $P\in \mathcal{X}\setminus \{ P_{\infty} \}$ such that $\ell$ contains the point $\varphi (P)=(a:c:1:c^{2})$ for some $a, c\in k$; namely, $\ell$ coincides with the line $\ell_{c}^{(2)}$.   

$\rm(I\hspace{-.01em}I)$ We consider non-skew Galois lines which are contained in the hyperplane defined by $Z=0$. We show that the line $\ell_{d}^{(3)}$ defined by $Y+dW=Z=0$ is such a non-skew Galois line for any $d\in \mathbb{F}_{q^{2}}$. If $d=0$, the field extension $k(\varphi (\mathcal{X}))/\pi_{\ell_{0}^{(3)}}^{\ast}k(\mathbb{P}^{1})$ coincides with the field extension $k(x, y)/k(y)$. Therefore, $\ell_{0}^{(3)}$ is a non-skew Galois line which is contained in the hyperplane defined by $Z=0$. We assume that $d\not= 0$. Let $G^{(d)}$ be the disjoint union of sets $A_{1}$, $A_{2}^{(d)}$, where
\[ A_{1}=\left\{ \left(
\begin{array}{cccc}
1 & 0 & \lambda & 0 \\
0 & 1 & 0 & 0 \\
0 & 0 & 1 & 0 \\
0 & 0 & 0 & 1
\end{array}
\right) \ \middle| \ \lambda \in \mathbb{F}_{q^{2}}, L_{q}(\lambda)=0 \right\} \]
and
\[ A_{2}^{(d)}=\left\{ \left(
\begin{array}{cccc}
1 & d^{-q} & a & d^{-2q} \\
0 & 1 & d^{-1} & 0 \\
0 & 0 & 1 & 0 \\
0 & 0 & d^{-2} & 1
\end{array}
\right) \ \middle| \ a\in \mathbb{F}_{q^{2}}, L_{q}(a)=(d^{-1})^{q+1} \right\} . \]
Let $M_{\lambda}\in A_{1}$ and $\sigma_{a}\in A_{2}^{(d)}$ be the matrices corresponding to $\lambda$ and $a$, respectively. Put $t=d^{-(q+1)}\in\mathbb{F}_q$. By direct computations, we have
\[ M_{\lambda}M_{\lambda'}=M_{\lambda+\lambda'},\quad M_{\lambda}\sigma_{a}=\sigma_{a}M_{\lambda}=\sigma_{a+\lambda},\quad \sigma_{a}\sigma_{a'}=M_{a+a'+t+t^2}. \]
Since $t^q=t$, we have $L_{q}(t+t^{2})=t+t^{q}=0$. Thus, the last matrix belongs to $A_1$. These identities imply that $G^{(d)}$ is an abelian subgroup of ${\rm Aut}(\mathcal{X})$. Since $|A_{1}|=|A_{2}^{(d)}|=q/2$, the order $|G^{(d)}|$ is equal to $q$. It is easily verified that the field extension $k(\varphi (\mathcal{X}))/\pi_{\ell_{d}^{(3)}}^{\ast}k(\mathbb{P}^{1})$ coincides with the field extension $k(x, y)/k(y+dy^{2})$. Since $\eta^{\ast} (y+dy^{2})=y+dy^{2}$ for any $\eta \in G^{(d)}$, it follows that $k(y+dy^{2})\subset k(\mathcal{X})^{G^{(d)}}$, where $k(\mathcal{X})^{G^{(d)}}$ denotes the subfield of $k(\mathcal{X})$ fixed by $G^{(d)}$. Since ${\rm deg}(\pi_{\ell_{d}^{(3)}})=q$, $k(y+dy^{2})=k(\mathcal{X})^{G^{(d)}}$. Therefore, $\ell_{d}^{(3)}$ is a non-skew Galois line which is contained in the hyperplane defined by $Z=0$. Conversely, let $\ell$ be a non-skew Galois line which is contained in the hyperplane defined by $Z=0$. By Lemma \ref{containinfinite}, $\varphi (P_{\infty})\in \ell$. Therefore, there exists $c, d\in k$ such that $\ell$ is defined by $cY+dW=Z=0$. If $c=0$, then $\ell$ is defined by $W=Z=0$. Then the field extension $k(\varphi (\mathcal{X}))/\pi_{\ell}^{\ast}k(\mathbb{P}^{1})$ coincides with the field extension $k(x, y)/k(y^{2})$, and $k(x, y)/k(y^{2})$ is an inseparable extension. This is a contradiction. Therefore, we can assume that $\ell$ is defined by $Y+dW=Z=0$ for some $d\in k$. If $d=0$, then $\ell$ coincides with the line $\ell_{0}^{(3)}$ defined by $Y=Z=0$. We assume that $d\not= 0$. Since the order of the hyperplane defined by $Y+dW=0$ at $P_{\infty}$ is equal to $1$, it follows that ${\rm deg}(\pi_{\ell})=|G_{\ell}|=q$. It is easily verified that the field extension $k(\varphi (\mathcal{X}))/\pi_{\ell}^{\ast}k(\mathbb{P}^{1})$ coincides with the field extension $k(x, y)/k(y+dy^{2})$. Since $|G_{\ell}|=q$, there exists $\eta \in G_{\ell}$ which is not contained in $A_{1}$. By Fact \ref{auto}, $\eta$ is represented by 
\[ \eta=\left(
\begin{array}{cccc}
1 & \zeta b^{q} & a & \zeta^{2}b^{2q} \\
0 & \zeta & b & 0 \\
0 & 0 & 1 & 0 \\
0 & 0 & b^{2} & \zeta^{2}
\end{array}
\right) , \]
where $a, b, \zeta \in \mathbb{F}_{q^{2}}, b^{q+1}=a+a^{2}+\cdots +a^{\frac{q}{4}}+a^{\frac{q}{2}}, \zeta^{q+1}=1$. Since $k(\mathcal{X})^{G_{\ell}}=k(y+dy^{2})$,
\[ (b+db^{2})+\zeta y+d\zeta^{2}y^{2}=(\zeta y+b)+d(\zeta y+b)^{2}=\eta^{\ast} (y+dy^{2})=y+dy^{2} \]
hold. Therefore, $b(1+db)=b+db^{2}=0$ and $\zeta =1$ hold. If $b=0$, then $\eta \in A_{1}$. This is a contradiction. Therefore, $b\not= 0$ and $d=b^{-1}\in \mathbb{F}_{q^{2}}$.

Finally, we describe the structure of the Galois group in case (3) when $d\not=0$. Put $t=d^{-(q+1)}$ and fix an element $\sigma=\sigma_{a}\in A_{2}^{(d)}$. By the multiplication formulas above, $\sigma^{2}=M_{t+t^{2}}$, $\sigma^{4}=1$. Since $G^{(d)}$ is abelian and $M_{\lambda}^{2}=1$, every element of $G^{(d)}$ has order at most four. Therefore, by the structure theorem for finite abelian groups, $G^{(d)}\cong (\mathbb{Z}/2\mathbb{Z})^{\oplus r}\oplus (\mathbb{Z}/4\mathbb{Z})^{\oplus s}$ for some non-negative integers $r$, $s$ satisfying $r+2s=n$. Since every element of $G^{(d)}$ belongs to $A_{1}$ or is of the form $M_{\lambda}\sigma$, we have $(G^{(d)})^{2}=\{1, \sigma^{2}\}$, where $(G^{(d)})^{2}$ denotes the subgroup of squares in $G^{(d)}$. On the other hand, the above decomposition implies that $|(G^{(d)})^{2}|=2^{s}$.

If $d^{q+1}=1$, then $t=1$ and $\sigma^{2}=1$. Hence $s=0$ and $r=n$. Therefore, $G^{(d)}\cong(\mathbb{Z}/2\mathbb{Z})^{\oplus n}$.

If $d^{q+1}\not=1$, then $t\not=1$. Since $t\not=0$, we have $t+t^2\not=0$, and hence $\sigma^2\neq 1$. Thus $s=1$ and $r=n-2$. Therefore, $G^{(d)}\cong (\mathbb{Z}/2\mathbb{Z})^{\oplus n-2} \oplus\mathbb{Z}/4\mathbb{Z}$.
\end{proof}

As an application of the main theorem, we determine all Galois points of the following plane curve.

\begin{coro}
The set of all Galois points for the curve defined by the affine equation
\[ y^{q+1}=x^{q}+x+(x+x^{2}+\cdots +x^{\frac{q}{2}})(y+y^{2}+\cdots +y^{\frac{q}{2}}) \]
coincides with the set $\{ (0:1:0), (1:1:0), (1:0:0)\}$. The associated Galois groups are $C_{q+1}$, $C_{q+1}$, $(\mathbb{Z}/2\mathbb{Z})^{\oplus n}$, respectively.
\end{coro}

\begin{proof}
We consider the projection $\pi_{P}:\varphi (\mathcal{X})\dashrightarrow \mathbb{P}^{2}$ from $P=(0:1:0:1)$. $\pi_{P}$ is represented by $(X:Y+W:Z)$. The field extension induced by $\pi_{P}$ coincides with $k(x, y)/k(x, y+y^{2})$. It is easily verified that $[k(x, y):k(x, y+y^{2})]\leq 2$. Note that $[k(x, y):k(x)]=q+1$. Since $q+1$ is odd, $k(x, y)=k(x, y+y^{2})$ holds. Therefore, $\pi_{P}$ is birational onto its image. Put $u=y+y^{2}$. Since $L_{q}(u)=y+y^{q}$, we have
\[
\begin{aligned}
u^{q+1}&=(y+y^{2})(y^{q}+y^{2q}) \\
&=y^{q+1}+y^{2q+2}+y^{q+1}(y+y^{q}) \\
&=L_{q}(x)+L_{q}(x)^{2}+L_{q}(x)L_{q}(u) \\
&=x^{q}+x+L_{q}(x)L_{q}(u).
\end{aligned}
\]
Therefore, $\pi_P(\varphi(\mathcal{X}))$ is contained in the plane curve defined by $u^{q+1}=x^q+x+L_q(x)L_q(u)$. Since the degree of the defining polynomial in $u$ is $q+1=[k(x, u):k(x)]$, the polynomial is irreducible over $k(x)$. Therefore, $\pi_{P}(\varphi (\mathcal{X}))$ coincides with the curve defined by the above equation. Since $\pi_{P}$ is birational onto its image, Galois points for the plane curve correspond to Galois lines for $\varphi(\mathcal{X})$ passing through $P$. If a line $\ell_{b}^{(1)}$ in the main theorem passes through $P$, then $(b+b^{2})^{q}=0$. Thus, $b=0$ or $b=1$ hold. A line $\ell_{c}^{(2)}$ in the main theorem does not pass through $P$. If a line $\ell_{d}^{(3)}$ in the main theorem passes through $P$, then $d=1$ holds. Therefore, there are exactly three Galois lines passing through $P$, defined by
\[ X=Z=0,\quad X+Y+W=Z=0,\quad Y+W=Z=0. \]
Their images under $\pi_P$ are $(0:1:0)$, $(1:1:0)$, $(1:0:0)$, respectively.
\end{proof}

\section*{{\bf Acknowledgments}}
The author would like to express his sincere gratitude to Professor Satoru Fukasawa for his invaluable advice, helpful discussions, and unwavering support and encouragement over the years.

\section*{{\bf Declaration of generative AI and AI-assisted technologies in the manuscript preparation process}}
During the preparation of this work, the author used ChatGPT (OpenAI) to improve the English language and readability of the manuscript. The author reviewed and edited the suggested revisions as necessary and takes full responsibility for the content of the manuscript.

\end{document}